\documentstyle{article} 

\newcommand{\ble}{\begin{lemma}}
\newcommand{\ele}{\end{lemma}}
\newcommand{\bth}{\begin{theorem}}
\newcommand{\eth}{\end{theorem}}

\newcommand{\ddx}{{\partial \over {\partial x}}}

\newcommand{\ddy}{{\partial \over {\partial y}}}
\newcommand{\ddz}{{\partial \over {\partial z}}}
\newcommand{\ddp}{{\partial \over {\partial p}}}

\newcommand{\ddth}{{\partial \over {\partial \theta}}}

\newcommand{\ddxbar}{{\partial \over {\partial \bar x}}}
 \newcommand{\ddybar}{{\partial \over { \partial \bar y}}}
\newcommand{\ddzbar}{{\partial \over {\partial \bar z}}}

\newcommand{\Ri}{Riemannian }
\newcommand{\R }{I \!\! R}

\newcommand{\bP }{I \! P}

\newcommand{\E }{{\cal E}}

\newcommand{\D}{ {\cal D} }

\newcommand{\De}{{\cal D ^{\epsilon}}}

\newcommand{\phiet}{{\phi^{\epsilon}_t}}
\newcommand{\phie}{{\phi^{\epsilon}}}

\newcommand{\Psie}{{\Psi^{\epsilon}}}

\newcommand{\pie}{{\pi ^{\epsilon}}}

\newcommand{\Ve}{{V_1 ^{\epsilon}}}

\newcommand{\Mth}{{M_{\theta}}}
\newcommand{\We}{{W^{\epsilon}}}

\newcommand{\LLe}{{{\cal L}^{\epsilon}}}
\newcommand{\elle}{{\ell^{\epsilon}}}

\newtheorem{proposition}{Proposition}
\newtheorem{theorem}{Theorem}
\newtheorem{backtheorem}{Background Theorem}
\newtheorem{lemma}{Lemma}
\newtheorem{definition}{Definition}

\newtheorem{question}{Question}
 
\title{\sc  Deformations of Nonholonomic Two-plane Fields
in Four 
Dimensions}

\author{ Richard Montgomery \\ Mathematics Dept. UCSC, Santa
Cruz, CA 95064 \\USA, email:  rmont@cats.ucsc.edu\\}
\begin{document}

\maketitle

{\bf Abstract:  
 An Engel structure is a maximally non-integrable
field of two-planes tangent to  a four-manifold. 
Any two Engel structures are locally diffeomorphic.
We investigate  the deformation space of
Engel structures obtained by deforming  
  certain canonical  Engel structures  on   four-manifolds
with boundary.
When the manifold is $RP^3 \times I$
where $I$ is a closed interval,  we show that this deformation 
space contains a subspace 
corresponding to Zoll metrics on the two-sphere (metrics
all of whose geodesics are closed) 
modulo 	`projective'
equivalence.   
The main tool is  a  construction of  an
 Engel manifold from a  three-dimensional contact manifold
  and a method of  reversing this 
construction.  These are special
instances of Cartan's method of prolongation and deprolongation.
The double prolongation of a surface 
X is an Engel manifold of the form $SX \times S^1$
where $SX$ denotes the unit tangent bundle to X.
The $RP^3 \times I$ example occurs
in this way, since the unit tangent bundle of
the two-sphere is $RP^3$.  
Besides proving these
new results,  the article has the flavour of a review.}

\section{Introduction}

The past three decades have seen the birth
and flowering of the 
field of contact (and symplectic)   topology
(\cite{Eliashberg}).
A key feature  giving contact structures
a  topological,  as opposed to purely 
 geometric,  nature is the  Darboux theorem.
This  asserts that any two contact (or symplectic)
structures of the same dimension are locally
diffeomorphic.  Hence contact structures have no local invariants.
 In the language of
singularity theory, contact structures are  {\sc stable
germs}.  (See the next \S, just preceding background theorem 2
for a definition of this term.)
Besides contact fields, their
even-dimensional analogues, and line
fields, the only other   stable germ  of 
k-plane fields in n-dimensions occurs
for 2-plane fields in 4-dimensions. 
The corresponding structures are called  Engel 
structures \cite{Vershik1}, \cite{BryantEtAl}. 
They are not nearly as well
investigated as contact structures.  

A basic theorem in contact  topology, referred
to as Gray's theorem, \cite{Gray},  
asserts that any deformation of
a global contact structure is diffeomorphic
to the original.  (See for example the
notes of Eliashberg \cite{Eliashberg} or
the text of Bryant \cite{Bryant}.)  
The analogous theorem is
false in Engel geometry,
as was  observed by Gershkovich
\cite{Gershk}.  The reason underlying
this failure is simple. One can canonically
associate to every Engel structure a certain line
field (see \S 4).
When the Engel planes are deformed so
is the line field. But line fields
on closed manifolds are well-known to
be topologically unstable.  Very recently,
Alex Golubev \cite{Golubev} has proved a version of Gray's
theorem in which the line field is fixed throughout
the deformation

The situation on manifolds with boundary is
more subtle. If
the Engel manifold is of the form $M \times I$,
$I$ an interval, 
with the line field tangent to $I$
then any small deformation of the line field
will yield a diffeomorphic line field.
  The purpose of
this paper is to investigate
deformations of certain canonical 
Engel manifolds of precisely this form. 
We prove that (A) any sufficiently
small deformation of this structure  
can be embedded into a larger structure
$M \times J$ of the same canonical form, but
despite this
fact, that  (B)   the  typical deformation is {\bf not}
diffeomorphic to the original.   We 
then relate the deformations
of a canonical structure on
$RP^3 \times I$ to
the deformations of the round metric
on the two-sphere through Zoll metrics.    

There are two surprises in this paper.  The first is
that in a semi-local sense the moduli space
of Engel deformations is equal to
the quotient space of the space of (germs of)
2nd order scalar ODEs 
$y'' = f(x,y,y')$ 
by (germs of) diffeomorphisms of the xy plane of dependent
and independent variables. 
The second surprise is that in some sense Engel geometry
is the product of three-dimensional contact geometry
with that of the real projective line. More 
precisely, an Engel structure induces
a transverse contact structure and
tangential real projective structure on
the four manifold.   The tangential structure 
is tangent to the canonical line field.   This 
second fact is   known to a few experts
 (Bryant and Hsu \cite{BryantHsu}, but 
we believe it is still worth emphasizing.  

{\sc Acknowledgements:}

I would like to acknowledge the
students of the geometry seminar
at Santa Cruz, especially Andrew
Klingler and Alex Golubev for
helpful comments. I would like to thank 
 Robert Bryant for discussions,  comments and
guides to the literature.  I would
like to thank Y. Eliashberg and M. Gromov
for raising versions of the questions posed
here. This work
was partially supported by NSF grant
DMS-9400515 and a Faculty Research Grant
from the University of California
at Santa Cruz.  

\section{Preliminaries.}

\subsection{Definitions}
We begin rather generally.  By a {\sc distribution}
we mean a smooth linear subbundle $\D \subset TQ$
of the tangent bundle $TQ$ of a manifold.
We will also call $\D$ is a k-plane field
where $k$ is the rank of the subbundle.
We may think of  $\D$   as a locally free
sheaf of smooth vector fields on the manifold.
Then we can use the notation 
$[\D, \D]$ for the sheaf generated by
all Lie brackets $[X,Y]$ of sections $X, Y$
of $\D$. Set $\D^2  =[\D, \D]$
and continue 
to take Lie brackets, setting $\D^{j+1} = \D^j + [\D,
\D^j]$. ( One can check that, as sheaves 
$D^j  \subset \D^{j+1}$.)

\begin{definition}
An {\sc Engel field} (or structure) on a four-manifold $Q$
is a rank 2 subbundle $\D \subset TQ$  of the tangent
bundle with the property that 
$\D^2$ is a rank 3 distribution and $\D^3$ is
the entire tangent bundle.
\end{definition}
Concretely, 
an Engel structure
on a four-manifold is a two-plane field
such that in a neighborhood of 
any point of the manifold  we can find a local frame $X,Y$
for the field such that $X,Y, [X,Y], [X,[X,Y]]$
is a frame for the tangent bundle.

\subsection{Prolongations}
 
Let $\E$ be a distribution of k-planes 
on a manifold $Q$.  The  projectivization 
$\bP \E$ is a fiber bundle over
Q whose typical fiber $\bP \E _q$ is
the real projective space of dimension
$k-1$ which consists of the one-dimensional subspaces
$\ell \subset \E_q$.   
The projectivization inherits a  
canonical distribution defined by declaring that a curve 
$(q(t) , \ell (t))$ in $\bP \E$ is tangent
to the distribution  
if and only if the derivative of the point of
contact, $q(t) \in Q$, lies in the line $\ell(t)$.
Alternatively, the distribution plane  at
the point $(q, \ell)$ is $d \pi _{(q, \ell)}  ^{-1} (\ell)$
where $\pi: \bP \E \to Q$ is the projection and
$d \pi$ its differential. We call $\bP \E$ with
this distribution the {\sc prolongation} of $Q, \E$.
This is a  special case of the
of Cartan's process of prolongation 
  \cite{Cartan}.   See final
section of 
 \cite{BryantHsu} ,
and other references to Cartan therein.

Instead of considering lines in $\E$ we could consider rays
in $\E$.
In this way we would obtain a
sphere bundle $S(\E)$ over Q,
with a distribution defined in the same way
as for $\bP \E$. Its distribution is the pull-back of 
of the distribution
just defined on $\bP \E$ by the 2:1 cover $S \E \to \bP \E$.
Either of the  spaces $\bP \E$
or $S \E$, together with
the  distribution  just defined
will be   called the {\sc prolongation of} $\E$.
When we need to differentiate between them,
we will call $S \E$ the {\sc oriented prolongation}. 
Observe that if $\E$ has rank 2, then the process of
prolonging increases the dimension of
the space by 1, but keeps the rank of the distribution
the same.

\subsection{Examples}

\subsubsection{ The basic Engel manifolds} 
 Let $M$ be a contact 
three-manifold with contact structure $\xi \subset TM$.
Then its  prolongation  $Q = \bP \xi$ 
and  oriented prolongation $S \xi$
form  the canonical
examples of  Engel manifolds.    We will
refer to $S\xi$ as {\sc the basic examples}.
We compute in \S 2.3.3 that these are indeed Engel
manifolds.  

We will be mostly concerned with the case
where $\xi$ is topologically trivial, which is
to say it admits two independent global sections.
In this case the prolongation $S\xi$
is equal to $Q \times S^1$ as a manifold.  

\subsubsection{The contact elements to a surface}

Take $X$ to be a two-dimensional surface
 and $\E$ to be its entire tangent space, $TX$.
Then its prolongation $\bP TX$ is called
{\sc  the space of contact
elements} to X.  It  forms a contact three-manifold
and is one of the most basic examples of  
a contact manifold.    
If X is oriented, then $STX$ can   be canonically identified
with $ST^*X$, the space of co-oriented contact elements.
A metric on
X defines a map between   the  oriented
prolongation $STX$ and the space of  unit tangent vectors to  X.
 
Combining this
with the previous example, we
see that   the double  prolongation of
a  surface X forms  an Engel manifold. 
This Engel manifold is 
diffeomorphic to $STX \times S^1$. 
 One of these sections is the  vertical field
relative to the fibration $STX \to X$.
Its flow rotates
tangent lines to the surface.  The other section depends on
choosing a metric on X.  Its integral   curves
correspond to the geodesics of this metric.

If we take $X = S^2$, the standard round two-sphere,
then $STX \cong SO(3) \cong RP^3$.  So
the double prolongation of the two-sphere yields
an  Engel structure on the Lie group $SO(3) \times S^1$.
It is invariant under left translations in this group.
The explicit identification $SO(3) \cong STS^2$
is given by
$g \in SO(3) \mapsto (ge_3, g e_1) \in STX$
where  $e_1, e_2, e_3$ are
the standard basis for $\R^3$.  
Let I, J, K be the standard basis for
the Lie algebra of the rotation group,
thought of as left-invariant vector fields.
K and I span   $\xi$. 
K is the vertical vector field 
relative to the projection $STS^2 \to S^2$.
The field I represents geodesic flow. 
(The field J represents a  different 
realization of   geodesic flow.   See \S 4.1,
especially the remarks following Theorem 4
and the review article of Arnold, \cite{Arnold2}
concerning these two realizations of geodesic flow.) A
global frame for the Engel structure is  
$\ddth$ and $cos (\theta) K + sin(\theta) I$
where $\theta$ is the angular coordinate on 
the $S^1$ factor.

\subsubsection{Coordinates}

According to the Darboux theorem,
centered at any point of a contact three-manifold
one can find coordinates $x,y,z$ such
that in this neighborhood the  contact distribution $\xi$
is described by the vanishing of the  one-form $dz - ydx$.
That is to  say:
$\xi_{(x,y,z)} = \{ dz - ydx  = 0 \}$.

Now $dx$ and $dy$ form {\bf linear}
coordinates on each contact plane
$\xi_{(x,y,z)}$.  Thus almost any   contact line
$\ell  \subset \xi_{x,y,z}$ is characterized 
by its slope, w relative to these coordinates.
(The only line not so covered is of course
the vertical line $dx = 0$.)
In other words:  $\ell = \{ dy = wdx \}$.
Thus $w$ forms an  affine coordinate
on the real projective line $\bP \xi_{(x,y,z)}$
and $x,y,z,w$ coordinatize the prolongation
$\bP \xi$.  
It follows from the definition in \S 2.3.1 that
the   Engel distribution on $\bP \xi$ (or $S \xi$ ) is
 the two-plane field  annihilated by the two one-forms $dz -ydx$
and $dy - wdx$.  One calculates that
the vector fields 
\begin{equation} 
W = {\partial \over {\partial w}}
\label{eq:fr1}
\end{equation} 
and
\begin{equation}
X = {\partial \over {\partial x}} +
 w{\partial \over {\partial y}} 
+ y{\partial \over {\partial z}}
\label{eq:fr2}
\end{equation}
frame this distribution.

\subsection{Background Results}

In this section we present  background
results on which we will be building,
and which are neccessary to understand
the importance of Engel distributions.  
 
 \begin{backtheorem}  [Engel normal form.]  
Any two Engel manifolds are locally
diffeomorphic.  More precisely, in a
neighborhood of any point of an Engel manifold
there exist coordinates $x,y,z,w$ centered at
that point such the distribution is spanned by
the vector fields W, X given by the equations
(\ref{eq:fr1}) and (\ref{eq:fr2})  above.
\end{backtheorem}

This theorem is  attributed to Engel.
It is proved several times 
in several different ways in the works of Cartan.  
See \cite{BryantHsu}, \S 2,  and references therein.
A modern proof of this theorem can be found 
in  the text  \cite{BryantEtAl}, Theorem II.5.1.   

 The conditions for a
distribution to be Engel are open conditions,
so that if we slightly perturb an Engel
distribution (relatively to the Whitney topology)
the resulting distribution is still Engel.  In particualar, it is locally
diffeomorphic to the original. 
A geometric  object is called {\sc locally stable} (in the
sense of singularity theory): if
any sufficiently small deformation of the object is
locally diffeomorphic to the original.
Thus Engel distributions are locally stable.
This is a property which they share with 
contact distributions.  
Unlike contact distributions, Engel distributions
are {\bf not} globally stable.  {\sc One of the
main aims of this paper is to understand and quantify
this failure of global stability.}

In order to explain the importance of Engel
distributions we  will   want another definition.
 \begin{definition}  A distribution is 
called {\sc regular} if for each
$j = 2,3, \ldots$ the dimensions of 
the spaces $\D^j (x)$ obtained by
evaluating all of the vector fields in
the sheaf $\D^j$ at the point $x$ are
constant. 
\end{definition} 

\begin{backtheorem}[stability theorem] 
Let $(k,n)$ be the rank of a distribution
and the dimension of the space in which it
lives.  Then the only stable distributions
occur when $(k,n) = (1,n), (n-1, n)$ or
$(2,4)$.  Any stable regular distribution
with  $(k,n) = (2,4)$ is an Engel distribution.  
\end{backtheorem}

A proof of this  can be found in Vershik and Gershkovich
\cite{Vershik1}, \cite{Vershik2}, and
also \cite{generic}. 

\begin{backtheorem}[Gray's theorem, \cite{Gray}]
Let $\xi = \xi_0$ be a contact distribution
and  $\xi_t$ a deformation of this
distribution through  contact distributions.
Then there is a one-parameter family of diffeomorphisms
$\phi_t$ of the underlying manifold such that
$\phi_t ^* \xi_t = \xi_0$
\label{backtheorem:gray}
\end{backtheorem}

Below, in lemma 10 of \S 9,
we give  a proof of a slight generalization
of this theorem.  As discussed above,
this theorem asserts that contact distributions
are globally stable.

Is the analogous theorem true for Engel manifolds?
NO! The reason is that   Engel manifolds inherit
a canonical line field, as we will now describe. 

\subsection{The Characteristic  line field} 
 
 An Engel distribution $\D$ determines a  
line field 
$$L \subset \D. $$ 
We will call it the {\sc characteristic line field}.
It is a central object in our investigations.
$L$ may be defined by
the   property  
$$[L, \D^2] \equiv 0  \; \hbox{ mod } \D^2 .$$  
In  canonical local coordinates
  $L$ is the span   of the vector field $W$ of
eq(1) above. In the BASIC EXAMPLES the integral
curves of   $L$   are
the fibers of the fibration $S\xi  \to M$. 
In other words, they are obtained by spinning the
contact line without moving the point of contact.
The corresponding vector field
will be written $\ddth$.

Because $L$ forms such a basic object it is worthwhile giving
another description for it which accentuates
its intrinsic aspect.  Given any distribution
$\D$, the Lie bracket $[X, Y]$
of vector fields X, Y
induces bilinear maps  $\D^j \times \D ^k \to \D^{j+k}$.
Set $V_j = \D^j / \D^{j-1}$.  
Observe that $[X, fY] = f[X,Y] \; (mod Y)$.  It
follows that the
Lie bracket induces  maps
$V_j \otimes V_k \to V_{j+k}$ of  a tensorial nature.
In particular if the distribution is regular so that
the $V_j$ are themselves vector bundles then
this induced bracket is a linear vector bundle map.
In the Engel case $V_2$ and $V_3$ are both one-dimensional
real vector bundles which we 
write (locally) as $\R$.  Thus
the Lie bracket induces   a bilinear
map $\D \otimes \R \to \R$.  Its 
one-dimensional kernel is $L \otimes \R$.
More concretely, if $X$ and $Y$ form a local frame
for $\D$ then $L_q \subset \D_q$ is the kernel of the map
$v \mapsto [[X,Y], v] (mod \D^2)$,
$v \in \D_q$.

{\sc remarks.}  The Engel line field enjoys several 
remarkable properties. 
Fix the endpoints of any
sufficiently short integral curve $C$ of $L$.
Consider the space of all unparameterized
 curves tangent to the
Engel distribution and joining these particular points,
and put the $C^1$-topology on this path space.
Then $C$ is an {\bf isolated} point
in this path space.  Bryant and Hsu 
 \cite{BryantHsu}  call this
property 
{\sc $C^1$-rigidity}.  

The second, related,  property is that the curve $C$
is a minimizing   subRiemannian (or Carnot-Caratheodory)
geodesic, 
independent of choice of metric on
the 2-plane fields $\D$.  See \cite{minimizers}, \cite{survey}, 
\cite{LiuSussmann}.  We will not have the occassion to
use these properties here though.

\subsubsection{Application: Deformations on Closed manifolds.}

Line fields are essentially vector fields,
An important, basic fact  in dynamical systems
that there are vast classes of vector fields which
are globally unstable: most nearby vector fields
are not conjugate to them.   Hence the following
result should come as no surprise:

\begin{backtheorem} [Gershkovich, \cite{Gershk}]
A typical deformation  of the
canonical Engel structure on $S \xi $ 
or $\bP (\xi)$ is not diffeomorphic
to the canonical structure.  
\end{backtheorem}

{\sc  proof:}
Suppose for simplicity that $\xi$ is parallelizable
so that $S \xi = Q \times S^1$.
 The characteristic line field associated to the  
canonical Engel structure on $S \xi $ is
spanned by the vertical vector field $\ddth$.
All the orbits of this vector field are periodic.
Fixing a value of $\theta \in S^1$ determines
a global section for this line field
and the corresponding Poincare return map is
the identity.  
This is not a generic property for line fields! 
When  
the Engel field is deformed, so is its canonical line
field.  This deformation is
rather arbitrary. In lemma 6 of \S 6 below
we show that {\bf any} contact map isotopic to the
identity can be realized as the Poincare return
map for some deformed Engel line field.  QED

\section{ Global Engel Deformations}

The last theorem asserts that Engel deformations
are not globally stable.
What can we save anything from the wreckage
of  global stability?
Since Engel distributions are
not globally stable, what is their
deformation space? 

{\sc question} [Courtesy of Viktor
Ginzburg.]  Suppose that the
 Engel  line field $L_t$ of a 
deformation $\D _t$ of Engel structures  remains constant.
Is $\D_t$ then diffeomorphic to $\D_0$?

This question was
answered very recently in the   affirmatively by Alexander
Golubev \cite{Golubev}.   As a special case
of his theorem we have: 
 \begin{backtheorem}[Golubev]
Let $\D_t$ be a   deformation of the canonical
Engel structure $\D_0$ on $Q = S \xi $
where $\xi$ is a {\bf parallelizable} contact structure
on a three-manifold.  Suppose that the
Engel line fields   $L_t$ are independent of t:
$L_t = L_0$.
Then there is a diffeomorphism
$\phi_t: Q \to Q$ taking 
$\D_t$ to $\D_0$.
\end{backtheorem}
 We could
provide a proof of Golubev's theorem in this
special case by slightly altering the lines of the
proof of theorem 1 below.  We choose not
to do this in the interests of space and refer the interested
reader to Golubev's preprint.

We now come to the main subject of this
article.  
Instead of fixing $L$, 
 we destroy its recurrence 
 by cutting out a
section of $Q = M \times S^1$ thus obtaining
a Engel four-manifold (with
boundary)  of the form $M \times I$,
$I$ an interval.   The line field
is tangent to the interval.  If we
perturb the field ($C^1$-) slightly we will obtain
a diffeomorphic line-field so that the answer stability
question is no longer obvious.  {\bf The main question
of the paper will be : to what extent does this `cutting'
restore the global stability of $\D$? } 
  
Henceforth we assume that $\xi \subset TM$
is a parallelizable contact structure.  
This is the same as an oriented contact structure
with a global non-vanishing section. 
We select  two Legendrian
direction fields 
$V_0, V_1  : M \to S \xi$
which are never collinear.  Together these define
a domain
$$\Omega = \Omega(V_0, V_1) \subset S\xi$$
consisting of all directions $V$ lying between them.
To be more specific, let us use the same symbol
$V$ for  a direction and for any vector 
whose positive span is this direction.  
Any direction $V$ can be expanded 
\begin{equation}
V= V(\theta, m) = cos \theta V_0 (m) + sin \theta V_1 (m).
\label{eq:T}
\end{equation}
uniquely  in terms of 
the two fields $V_0, V_1$.  
The domain  $\Omega $ consists of those directions $V$ for
which $0 \le \theta \le {\pi \over 2}$.

\begin{definition}  The domain 
$\Omega = \Omega(V_0, V_1)$ just defined
will be called the {\sc standard Engel domain}
determined by the pair  $V_0, V_1$
of Legendrian direction fields.  
\end{definition}  

The expression 
(\ref{eq:T}) 
defines a global trivialization
$S \xi \cong M \times S^1$.
This trivialization maps
$\Omega$ diffeomorphically onto $M \times [ 0 , {\pi \over 2} ]$.  
The image of the
section $V_0$ corresponds to $M \times 0$ and
will be called the {\sc bottom} of the domain.
Similarly, the image of $V_1$ will be called the
{\sc top}.  The maps $m \to V(\theta, m)$ defined
by (\ref{eq:T}) are Legendrian direction fields
$$V(\theta) : M \to S \xi$$
 interpolating between $V_0$ and $V_1$.
Their   images will be denoted
$$M_{\theta} = image(V(\theta))$$
and correspond to $M \times \{\theta \}s$.

By using our trivialization,
we may also think of 
expression (\ref{eq:T})
as defining a vector field
(or direction field) on
$\Omega$ which is tangent to   the foliation by the $M_{\theta}$s.
Let   $\ddth$ denote the vertical
vector field on $S\xi$ which is tangent to 
 the $S^1$-factor of 
$M \times S^1$ under our trivialization.
Then $V$ together with $\ddth$ span the
Engel field. The Engel line field $L$ is
spanned by $\ddth$t.

\begin{question}  Suppose that 
we deform the Engel distribution
$\D = \D_0$ on a standard Engel domain
$\Omega$ into a new Engel distribution
$\D_t$.  Is $(\Omega, \D_t)$ diffeomorphic
to $(\Omega, \D_0)$?.  If not,
how can one describe the space of
nontrivial deformations?
\end{question}

{\sc \bf Answer:}
    
\begin{theorem}  
For all sufficiently small t
there is an Engel isomorphism
$\psi_t : (\Omega , \D_t) \to (\Omega_t, \D_0)$ 
between the standard Engel
domain  with varying Engel structure $\D_t$
and a varying Engel domain $\Omega_t = \Omega(V_{0 t}, V_{1 t})$
with standard Engel structure $\D_0$ induced
from the inclusion of $\Omega_t$ into $S \xi$.
In other words, we can ``straighten out'' the Engel deformation
at expense of varying the  top and bottom of the domain.
 
If, during the course of deformation,  the structure is  constant
in some neighborhoods of the top and bottom,
then we may take $V_{0t} = V_0$ and 
$V_{1t } = \phi_t ^* V_1$ for some contact
diffeomorphism $\phi_t: (M, \xi) \to (M, \xi)$.
In other words, the bottom remains unchanged and the top
varies by a contact transformation.  This
transformation $\phi_t$ is canonically defined by the
deformation and will be 
called the {\sc bottom-to-top map}.  
\end{theorem}

\begin{theorem}  Let
$\Omega_0 = \Omega (V_0, V_1)$ and
 $\Omega_1 = \Omega (W_0, W_1)$ be
two standard Engel domains in $S \xi  \cong M \times S^1$.
Then there is an orientation preserving  Engel diffeomorphism 
$\Omega_0 \to \Omega_1$ if and only
if there is a contact diffeomorphism of 
$M, \xi $ which takes the pair $(V_0, V_1)$ of 
Legendrian direction fields
to the pair $(W_0 , W_1)$.
 \end{theorem}

{\sc remark} When we say that
the Engel diffeomorphism preserves orientation
we mean that it preserves the orientation of
both the manifold and the distribution.
  See \S 3.4 for more on this.

\begin{lemma}[Realization lemma]
Let $\phi_t: M \to M$ be any 
sufficiently $C^1$-small contact isotopy
and $V_0, V_1$ any pair of everywhere independent
Legendrian vector fields for $(M, \xi)$.
Then there is a deformation $\D_t$ of the canonical Engel
structure on $\Omega(V_0, V_1) \subset S \xi $ for which
$\D_t = \D_0$ near the top and bottom of the domain
and such that the
bottom-to-top map induced by the deformation is 
$\phi_t$.  
\end{lemma} 

{\sc Remark}  ``Sufficiently small'' in the realization lemma
can be replaced by a certain `positivity' condition
concerning the sense of twisting of $\phi_t$.
This positive twist property may 
prove to be important in studying questions
regarding extensions and existence of Engel structures.  

Combining the lemma with the previous two
theorems, we have a  complete
translation of our
Engel deformation problem into contact terms.
Namely, the problem of  understanding  the 
Engel deformations of the Engel domain $\Omega(V_0, V_1)$
which are trivial near the top and bottom 
is equivalent to the problem of understanding
certain deformations $(V_{0t}, V_{1t})$ of
the original Legendrian line field
pair $(V_0, V_1)$ on $(M, \xi)$.
The  deformations are restricted to be of the form
$V_{0t} = V_0$ and   $V_{1t} = \phi_t ^* V_1 $ 
for some contact isotopy $\phi_t$.  
The  deformation is considered trivial if
we can find a contact isotopy $\psi_t$
which preserves $V_0$ (i.e. 
such that $\psi_t ^* V_0 = V_0$) and such that 
$\psi_t ^* V_{1t} = V_1$.  Two
such deformations $(V_0, V_{1t})$
and $(V_0, \tilde V_{1t})$ are equivalent
if we can find a contact isotopy preserving
$V_0$ and taking $V_{1t}$ to $\tilde V_{1t}$.

As we will see in the next section, such contact
isotopies are very rare.  Generically
two Engel domians are not diffeomorphic
and the moduli of such domains involved
functional parameters.    

\section{The Appearance of  Differential Equations}

\subsection{General 2nd Order Differential Equations}

We have reduced our
deformation
question to a question in contact
geometry concerning a pair of Legendrian
direction fields.  
Our  method of reduction
was formalized by Cartan and
is referred to as deprolongation in \cite{BryantHsu}.
We will now see how a further deprolongation
reduces this contact question
to a question regarding 
2nd order differential equations on   surfaces.

Consider a pair $(V_0, V_1)$ of direction fields
on a three-manifold $M^3$.  Construct the 
local quotient of $M^3$ by the integral curves
of $V_0$.  By ``local quotient'' we mean that
we restrict the field $V_0$ to a small enough neighborhood
$U$ so that the quotient forms a smooth
two-manifold $X^2$ with the quotient map
$\pi: U \to X^2$ a submersion.  For example, U could
be a flow-box for $V_0$.  Now use the
projection $\pi$ to push the 
integral curves of $V_1$ down to $X^2$.
The result is a one-parameter family of curves passing
through each point of $X^2$.  
If the pair spans a contact field 
in $M$ then the 
curves in X passing through a given point
will (microlocally) be in
one-to-one correspondence with tangent directions
through that point.  In other words, as we
move along an integral curve for $V_0$,
the tangents of the projections of the $V_1$-curves
also move, to first order.  A family of
curves parameterized by their
tangents is nothing more than a 2nd order differential
equation.   

The discussion above  paraphrases that of Arnold, \cite{Arnold},
(especially p.  52-53) and  Cartan, \cite{Cartan2},
(beginning on p. 25).  We find it helpful  
to make this discussion more explicit.   
\begin{lemma} [Normal form theorem]  
Suppose we are given a  of 
direction fields on  a three-manifold 
whose span is a contact field.  Then we can
find,  in a neighborhood of
any point, local coordinates $x,y,p$
and  vector fields  $V_0$, $V_1 $ spanning the pair of
direction fields such that: 
$$V_0 = \ddp$$
$$V_1 = \ddx + p \ddy + f(x,y,p) \ddp.$$ 
where $f$ is a function of $x,y,z$ 
which can be chosen so
as to vanish  at the origin $(0,0,0) = 0$.
We call this the {\sc normal form} for the pair.
The contact form annihilating
$V_0$ and $V_1$ is $dy - p dx$. 
\end{lemma}

Our appendix contains a
proof of this lemma.
  It follows from the lemma
that any regular Legendrian curve transverse to $V_0$
can be parameterized as 
$x \mapsto (x, y(x), p(x))$
where $p(x)  = {dy \over dx}$.   Now we proceed
with the construction above.  Project the 
integral curves of the {\bf second line field} $V_1$
 along the leaves of $V_0$ 
onto this quotient space.    
According to the   normal
form for $V_1$, the resulting curves must satisfy the
2nd order differential equation
\begin{equation}
{{d^2 y} \over {dx^2}} = f(x, y,p).
\label{eq:star}
\end{equation}

To reverse this procedure, suppose we are 
given a system of 2nd order differential equations
 on a two-dimensional surface $X$.
We are  interested in ``parameterization independent''
differential equations.
The Euler-Lagrange equations for
a Lagrangian $L: TX \to \R$,
 homogeneous of degree one
in the velocities, provide a class of example.
This class of examples contains 
the  geodesic equations for  metrics on $X$.  
In modern terms, 2nd order differential equations 
are usually described in terms of a  {\bf spray}.
(Alternatively, they could be described as sections  of 
the bundle $J^2 X \to TX$
where $J^2 X$ denotes the 2nd jet space for
maps of the real line into X.  Recall that
the tangent space $TX$ is the 1st jet space for
such maps.)  
In  local coordinates $(x,y)$ 
such a  system of second-order differential equations 
is written
 ${d^2  x \over dt^2} = f_1( x,y, \dot x, \dot y)$,
${d^2  y \over dt^2} = f_2 ( x,y, \dot x, \dot y)$.
The parameterization-independence condition
implies that 
$f_1 (x, y, \lambda \dot x, \lambda \dot y) = 
\lambda^2 f_1( x,y, \dot x, \dot y)$ and that 
$f_2 (x, y, \lambda \dot x, \lambda \dot y) = 
\lambda^2 f_2( x,y, \dot x, \dot y)$ for 
$\lambda$ a positive number.  
However these are not sufficient conditions.
(For example, all the solutions
to a standard form for the geodesic
equations, the form in which the
Lagrangian is the kinetic energy,  are  parameterized by a
constant multiple of arc-length, and so the equations
are not parameter-independent.)   

A more concise way to write down {\bf parametrization-independent}
2nd order differential equations on the surface is 
the one described above.  We  take $M = STX$, $\xi$ to be the
standard contact form, and $V_0$ to be the vertical direction field.
(See \S 2.3.2.)
Then a  Legendrian direction field $V_1$ which is transverse to $V_0$
describes a 2nd order differential equation on $X$
whose integral curves have, a priori,  
no distinguished parameterization.  The Legendrian condition
insures that the integral curves of $V_1$ are ``holonomic''
-- they are tangent to the appropriate contact elements.
The fact that $V_1$ is a direction field and is defined
on $STX$ rather than $TX$ implies that the integral
curves have no distinguished parameterization.
(However they do have a distinguished orientation.
If we wanted unoriented curves then we would
do better to work with line fields on $\bP TX$.)
The fact that $V_1$ is transverse to $V_0$
implies that for each tangent direction
there is a corresponding integral curve.  
 
Pick a particular contact element
$\ell \subset T_p X$ to $X$.
Choose coordinates $(x,y)$ centered at
the point $p$ and such that the   x-axis 
in these coordinates  represents
this contact element $\ell$.  Any $C^1$ curve on $X$
whose initial contact element
is close to $\ell$
can then be parameterized, near $p$,  in the
form of a graph:  $y = y(x)$.
If we use this
parameterization 
then the system of parameterization-independent
differential equations becomes
a single 2nd order scalar differential equation
as above, with $p = dy/dx$.   

We  call two parameter-independent differential equations
on a surface {\sc projectively equivalent}
if there is a diffeomorphism of the surface which
takes the solutions of one {\bf thought of as
unparameterized curves} to the solutions of the other.

{\sc \bf Example:}
Central projection is a local projective equivalence
between the geodesics on the round sphere and the 
geodesics on the flat plane.  

Projective equivalence of differential equations
boils down to contact equivalence between
Legendrian line-field pairs.  
We summarize the above  discussion as a background theorem.  
\begin{backtheorem} A  pair of 
 direction fields $(V_0, 
V_1)$ on a three-manifold  whose span
is a contact field     
locally defines a parameter-independent 2nd order 
differential equation on a surface.    And associated
to every 
such  equation
is a pair of Legendrian direction fields.   The  two
pairs of direction fields are locally diffeomorphic 
if and only if their   corresponding equations are projectively 
equivalent.  
\end{backtheorem}

We have thus reduced our 
Engel deformation question to one
regarding deformations of 2nd order differential
equations, a theory well-developed near the 
turn of the last century.  
( See Arnol'd (loc. cit.) and Cartan (loc.cit.).)
The crucial fact for us is that the space of
germs of parameter-independent 2nd order differential equations
modulo projective equivalence is infinite-dimensional.

Consider now  an Engel deformation
$\D_t$ of the standard Engel structure
$\D_0$ on the domain $\Omega(V_0, V_1)$
whose support is some small neighborhood $U$
of a point
$(m, \theta)$ relative to our trivialization.
It follows from the proof  of  theorem 1,
that the direction fields $V_{0t}, V_{1t}$
described there  
equal  $V_0, V_1$ except in the small neighborhood,
  $\pi(U)$ of the point $m \in M^3$.
This deformation
of fields in turn defines
a deformation of 
  parameterization-independent second order differential
equations, with the support
of the deformation being the quotient
of $\pi(U)$ by the leaves of $V_0$.

\begin{theorem}  
The bottom $V_0$, and top $V_1$
of the standard Engel domain
$\Omega(V_0, V_1)$, together
with a point $m \in  M^3$
define a germ of a parameter-independent
second order differential
equation which can be put into the  form (\ref{eq:star}).
Consider a point $q \in \Omega$ with $\pi(q) = m$,
and the corresponding space of Engel deformations
with support 
in a small neighborhood of $q$.  The germ (at $q$)
of this space of Engel deformations of $\Omega(V_0, V_1)$
  modulo Engel   isotopies 
is  isomorphic to  the space of  germs of 
deformations of this 2nd order differential equation
modulo projective equivalence.  
This latter space is infinite-dimensional.  
\end{theorem}

{\sc remark}  As can be
seen from an inspection of
the proof of theorem 1, 
this theorem is also true if we consider
an apparently larger class of  germs --
namely germs 
along an integral curve of the Engel line field.
In other words, consider deformations supported in
arbitrarily small neighborhoods of the form
$\{m \} \times I$ in our trivialization,
and impose the usual equivalence relation
on such deformations in order to
obtain germs of deformations along a leaf.

The point of the next
subsection is to give a global version
of these musings on differential equations.

\subsection{Geodesics and Zoll Surfaces.}

We refer the reader back to \S 2.3.2. 
Take $M = STX$ to be the space of oriented contact
elements for an oriented \Ri surface X with metric g.  
Take $V_0$ to be the vertical direction field for
the fibration $\pi: STX \to X$.  
Take $V_1$ to be the direction field corresponding to 
the  geodesics on X. The  projections of 
the  integral curves of $V_1$ along the
$V_0$ curves  are  geodesics on X
for the metric $g$.  The  geodesics
are to be thought of as
  unparameterized oriented curves.  In this manner the Engel domain
$\Omega(V_0, V_1)$ encodes the `unparameterized geodesic flow'
for X.    

Recall that a {\sc Zoll metric}
on the sphere is a metric all of whose geodesics are closed.
Two \Ri metrics are said to be {\sc projectively equivalent} if
there is a diffeomorphism which takes the geodesics of one
to the geodesics of the other, {\bf where the geodesics
are to be viewed as unparameterized curves}.  For example,
as discussed earlier, 
central projection from the
hemisphere to the plane establishes a local projective equivalence
between the round sphere and the plane.

\begin{theorem}   Let
$g_t$ be a deformation of the
round metric $g_0$ on the sphere through Zoll metrics. 
Let  $V_{1t}$  be the  corresponding
deformation of Legendrian direction fields on
$STS^2$ which  
describe the (unparameterized)  geodesics  for $g_t$.
Then
\begin{itemize}
\item{}(A)  $\Omega(V_0, V_t)$ 
can be realized as an Engel deformation
of the  standard domain $\Omega(K,I) \subset
SO(3) \times S^1\cong STS^2 \times S^1$
corresponding to a  deformation which is trival near the top
and bottom,  
\item{}(B)
the Engel deformation of (A)
is trivial if and only if
the metrics $g_t$ are isotopic 
to the round metric on a sphere $S^2 (r_t)$
of radius $r_t$, 
and 
\item{}(C) two such `Zoll' deformations are equivalent
if and only if the corresponding metrics $g_t$ and
$\tilde g_t$ are projectively equivalent.
\end{itemize}
\end{theorem}

{\sc Discussion.}  
Theorem 4 asserts that the Engel deformation space
for $\Omega(K,I)$ contains
a copy of the ``Zoll deformation space''
consisting of (Zoll metrics)/(projective equivalence).
To get a rough idea of this space,
first consider the perhaps more well-known
space of   Zoll metrics modulo
isometric equivalence which was investigated 
in the last few decades by Guillemin and others.
(See 
\cite{Guillemin}, \cite{Weinstein}
and \cite{Besse}.)  It is infinite-dimensional
and Guillemin
showed that its tangent space at the round metric consists
of `one-quarter' of the functions on the sphere. 
Projective equivalence is a coarser equivalence relation
then isometric equivalence, so 
that we have a  map (Zoll metrics)/(isometric
equivalence)  $\to$ (Zoll metrics)/(projective 
equivalence). The fiber of this map consists of the space of all  
 metrics  projectively equivalent to a given Zoll metric,
modulo the equivalence relation induced by isometry.
These fibers were investigated by Beltrami.  
(See the discussion in \cite{Darboux}.) 
They are discrete unless the given
metric has constant curvature in which case
the fiber is one-dimensional with the single one-dimensional parameter
corresponding to homotheties of the metric (or, what is the same, 
to the constant curvature). 
It follows that the space of Zoll metrics
modulo projective equivalence is also an infinite dimensional
space.   

The real
content of  theorem 4
is assertion (A).  Assertion (B)
follows immediately from (A) and the
discussion of the previous section
concerning differential equations.
Assertion (C) is classical and can be
found in the book of Darboux \cite{Darboux}.

In view of the realization lemma,
 (A) asserts that  
 the geodesic flows $V_t$
for  Zoll
deformations can be written in the
form $\phi_t ^* V_1$ for some contact
isotopy $\phi_t$.  This should be compared to
a   result of Weinstein.  (\cite{Weinstein}, esp. the final remark
and the proof of theorem 5.1.  See also \cite{Guillemin},
Appendix B.) 
Weinstein's theorem makes this same assertion,
but for a different representation of
geodesic flow.   
Weinstein
used the Hamiltonian framework for geodesic flow.
There the geodesic flow is defined as
the   Hamiltonian vector field for the Hamiltonian
corresponding to kinetic energy.
Any nonzero level set of this Hamiltonian  
 can be identified with $ST^* X$.  The resulting
vector field is  a Reeb vector field.  It is
transverse to the contact field.  
In contrast,
for our theorem the geodesic `flow' consists of
the integral curve of a {\bf Legendrian} direction field.
There can be no contact automorphism taking
the representation of Weinstein to the
representation which we use (Darboux's representation.)
Geometrically, the difference between the two is
that in Weinstein's representation  
momenta represent tangent elements to wavefronts and 
as such are perpindicular
to the direction  of the geodesic's motion,
whereas in our (Darboux) representation
momenta represent tangent elements to
the geodesics themselves.

\section{ Contactifying an Engel manifold} 

We  constructed
our basic examples of Engel
manifolds (\S 2.3.1 ) from three-dimensional contact manifolds.
Here we  reverse this  
procedure by dividing out
by the orbits of the Engel line field.
This is an instance of `deprolongation'
in Cartan's sense.    See the example
titled `Prolongation and Deprolongation'
in the final section of 
 \cite{BryantHsu} ,
and other references to Cartan therein.      

\subsection{Contactification} 
 
Let $(Q, \D)$ be an Engel manifold
and let $\cal L$ denote the foliation of $Q$ by
leaves of the canonical line field $L \subset \D$.   
Let  $M= Q/ {\cal L}$ be the quotient space. Its
points are the leaves of the foliation.
We will assume that it 
is a {\sc nice} quotient, by which we  mean that it 
is a  manifold and that the quotient map
$\pi: Q \to M$ is a submersion.   

The rank 3 distribution $\D^2 = [\D, \D]$ is invariant under the flow
along $\cal L$.  This is because the line field 
$L$ is
defined by the condition $[L, \D^2] \subset \D^2$.
We can thus push $\D^2$ down to the quotient
otainining a distribution $\xi = \pi_* \D^2$.
This quotient distribution
has rank 2 because  $L$ is the kernel of the
differential $d \pi$.

\begin{lemma} [Contactification lemma]  If the quotient $M = Q/
{\cal L}$ is nice (as defined above) then it is a contact manifold
with contact structure $\xi = \pi_* \D^2$.
We call it the {\sc contactification} of $Q$.
\end{lemma}

{\sc Proof:}  It only remains to prove
that $\xi$ is contact.  This follows immediately
from the fact that $[\D^2 , \D^2]$ is the
entire tangent space to Q. Alternatively,
we can use  the local normal form
(1) , (2). In terms of those canonical coordinates
$x,y,z, w$, the leaves are the $w$-lines
so that $x,y,z$ coordinatize the
quotient $M$.  The Engel distribution is defined locally
by the Pfaffian system
$\theta_1 = \theta_2 = 0$ where
 $\theta_1 = dz - ydx$ and  $\theta_2 = dy -wdx$.
(See also the proof of the contactification theorem
above.)
The form $\theta_1$ is the annihilator of
 $\D^2$  and pushes down
to a contact form on $M$. 
 
{\sc Example: }  Return to the basic
example.  The original
contact manifold $M$ 
can be identified with the quotient of $Q$
by the one-dimensional foliation.
The contact structure induced by
this construction is
the original contact structure.

\subsubsection{Transverse Contact Structure}

Let   $S \subset Q$ be a three-dimensional
submanifold of our Engel manifold  
transverse to the line field $L$.  We
think of $S$ as a small disc.  
Then $S$ may be identified with the
quotient space $U/ { \cal L}$
where $U$ is some (perhaps small)
neighborhood of $S$.  For
example, pick 
a vector field spanning   $L$ and then using the
flow of this vector field to 
 sweep out $U$.   
S inherits a contact structure by
identifying it with the quotient. This
structure is
simply the intersection of $TS$
with $\D^2$.  

Let $\tilde S$ be another contact slice,
and suppose that there is a leaf
$\ell$ of $\cal L$ which intersects
both $S$ and $\tilde S$. Choose a tube $U$
around $\ell$, that is to say,  
neighborhood $U$ of $\ell$ consisting
of arcs of leaves.
Now $U \cap S$ and
$U \cap \tilde S$
are both realizations of the
 the contactification
of $U$.  The flow along
the leaves in $U$ induces a contactomorphism
between these two realizations.  
This situation is summarized by saying that
the foliation $\cal L$ has a {\sc transverse contact
structure} This means that 
(A) any local
transverse to a given leaf inherits a contact structure,
and 
(B) these contact structures are consistent:
 the flow
along $L$ defines a map between   two
 different  local transverses to the same leaf 
which is a  contact map.

\subsection{Development Map}

Let $(M, \xi)$ be the  contactification
of an Engel manifold $Q$, as in \S5.1.1,
so that $M = Q/ { \cal L}$. 
Although the distribution $\D ^2$
descends to a  distribution on $M$ (the contact
distribution) the Engel field
$\D$ itself does not descend.  Define the map 
$$\Psi: Q \to \bP \xi $$  
  by mapping the point $q \in Q$
to the image of $\D_q $ in $\xi_{\pi(q)} = \pi_* \D^2 _q$.
 In other words
$$\Psi(q) = d \pi (q) (\D (q)) \subset \xi (\pi (q)).$$
Note that the image of $\D (q)$ under $d \pi$ is indeed
a one-dimensional subspace of $\xi (m), m = \pi (q)$
because the kernel of $d \pi$ is the line field $L$.
We call $\Psi$ the {\sc development map} for Q.

Restricting $\Psi$  to a leaf $\ell = \pi ^{-1} (m) \subset
Q$, of the line field $L$ we obtain a map
$$\ell \to \bP^1 = \bP (\xi (m)).$$ 
This map is monotone, i.e. its derivative is
everywhere nonzero. This is merely the
condition that that $[L , \D] \ne \D$, which
must alway hold, since $\D^2$ has rank 3.  

This construction gives each leaf $\ell$
a {\sc real projective structure}.  Let w denote
an affine coordinate on the line
$\bP^1$.  We may pull w back to the leaf
thus obtaining a coordinate on part of the leaf.
If $\bar w$ is another affine coordinate,
then it is related to $w$ by
a linear fractional transformation:
$\bar w = (a w + b) /(cw + d)$.  Each
leaf can be covered by such coordinates.
As we vary from leaf to leaf the coefficients
$a,b,c,d$ change but in such a way that they are
functions of coordinates $x,y,z$
on the leaf space M. 

The fundamental theorem
of projective geometry states that a projective
map from one projective line to another is determined
by its values on three points.  It follows
that if $\Psi:Q \to Q$ is an Engel isomorphism,
then the values of $\Psi$ along any leaf $\ell$
are completely determined by its values on
any small interval of this leaf.    

Let us summarize this discussion so far
\begin{lemma}
Let $(M, \xi)$ be the contactification of an Engel manifold
$Q$.  Then there is a canonically defined Engel immersion,
$$\Psi: Q \to \bP \xi,$$ which we call the
{\sc development} of $Q$. If $Q$ is closed then
the development is a covering map.  
In any case, it maps the leaves of the canonical line
field on $Q$ to the fibers of
$\bP \xi \to M$. 

If $\phi$ and $\psi$  are two Engel immersions 
whose values  agree   at three distinct but sufficiently
close points of a leaf of the canonical line
field, then they are  
equal everywhere along this leaf. 
\end{lemma} 

\subsubsection{Tangential  Projective structure}
The real projective structure on 
the leaves of $L$ is well-defined   whether or not the global quotient $Q/{ \cal L}$
is nice.  To construct it,   use local slices to the leaf in
order to define a local development map.  The differential
of the  contact map
intertwining the two slices defines a linear fractional
map between the projectivizations
of the contact planes at the leaf, and consequently
a linear fractional transformation between the corresponding
projective coordinates.  Consequently the leaf
has an atlas whose overlaps are linear fractional
transformations.   We summarize this
by saying that $Q$ has a {\sc tangential projective structure}.

\subsubsection{Oriented Development}

If $\D$ and $Q$ are oriented
then so are $L$ and $\D^2$. (See \cite{Gershk} or
\cite{meAnd}.) 
It follows
that $M$ and $\xi$ inherit natural orientations.
In this case it makes sense to formulate
an oriented version of the above proposition.

\begin{proposition}
If $M$ is the contactification of an oriented Engel manifold
$Q$ (assumed  ``nice'' as before) 
then there is a  canonically defined 
oriented Engel immersion
$$\Psi: Q \to S \xi ,$$
which we will again refer to as the {\sc development}
or {\sc oriented development} of $Q$.   
\end{proposition}

\subsubsection{Orientations}

Suppose that the manifold $Q$ and the Engel field $\D$ are
both oriented. The 
action of Lie bracketing  induces
orientations in both $L$ and $\D^2$.
 See \cite{Gershk} or \cite{meAnd} for
 details of this orientation construction.
If $M$ is an oriented three
manifold with oriented contact structure
$\xi$  then   $S \xi $ and its
Engel distribution are both  
canonically oriented.  Thus the
line field is canonically oriented. 
Under this  orientation of the characteristic line
field positive motion corresponds to
  rotating the contact directions
in the counterclockwise sense.  
It follows that any Engel diffeomorphism
preserving  orientations and
mapping one standard domain
$\Omega(V_0, V_1)$ to another,
$\Omega(W_0, W_1)$ must take  bottom ($V_0$)
 to  bottom ($W_0$) and   top ($V_1$) to 
top ($W_1$). See theorem 1 above
and the statement following it, where this
was alluded to.

\subsection{Contact slices}

\begin{definition}  A {\sc local cross-section}
or {\sc slice}
to an Engel manifold is an embedded three-dimensional
sub-manifold which is transverse to the Engel
line-field.  A {\sc global cross-section} 
or {\sc global slice}
is a local cross-section which intersects every leaf
of the line-field.
\end{definition}

{\sc Example:}  In the oriented basic example,
$Q = S \xi $, 
the submanifolds $\Mth$ are global slices.

\begin{lemma}.  Let   $M \subset Q$ be a cross-section
of the Engel manifold $Q$.  Then $\xi = \D^2 \cap TM$
defines a contact field on $M$.   The projectivization 
$\bP \xi$ of
this contact field admits a global section over $M$,
namely $Z = \D \cap TM$.  Any sufficiently small neighborhood of  
$M$ in $Q$ is diffeomorphic as an Engel manifold to a 
neighborhood of $Z \subset \bP \xi$ with its canonical Engel
structure.
\end{lemma}

{\sc proof:}  The slice condition 
implies that there is a  neighborhood of $U$ of
$M$ such that $U/{{\cal L}_U} \cong M$. 
The development map induces the Engel isomorphism
of neighborhoods.  QED 

Suppose now that the line field $L$ is oriented,
  that 
$Q$ is a closed compact manifold,  and that
$M$ is a global slice.
Define a Poincare return map
$ M \to M$ by mapping $m \in M$ 
to the next intersection point of the
leaf  through $m$ with $M$.  
We are guaranteed that such a point exists
because of the compactness of $M$ and $Q$ 
and the assumption that $M$ is a global slice.
This map can be realized as the flow of a vector
field along $\L$ and this flow preserves $\D^2$.
It follows from this and the lemma 2 
that the Poincare return map is a contact map
of $M$.

As a variation on this theme, suppose that  
 $M$ and 
$M ^{\prime}$ are two connected global contact slices
which do not intersect, and suppose that  the 
 characteristic line field
is oriented. 
Flowing along the leaves of $L$ now defines a contact  map
$M \to M ^{\prime}$.  
 The model example here
is the
domain $\Omega(V_0, V_1) \subset S \xi $,
with $M$ and $M^{\prime}$ being $M_0$ and $M_{\pi/2}$.

\section{Proofs}

\subsection{ Proof of the first half of theorem 1.}

  We apply the constructions of the last
section to prove the assertions of the first paragraph of   
theorem 1.  
  
Recall that the  
$M_{\theta} = image( V(\theta)) $
are global slices for the unperturbed Engel structure. 
In terms of the isomorphism of $S \xi $ to $M \times S^1$
these slices are obtained by setting the $\theta$ coordinate
equal to a constant.  As contact manifolds, they are all
isomorphic to $M$ with the section $V(\theta)$ providing 
the
isomorphism.

Now consider a deformation $\De$ of the canonical Engel
structure having the property that
$\De (q) = \D (q)$ for $q$ near the ``top'', $M_{\pi/2}$
and the ``bottom'', $M_0$ of $\Omega$.  For small
$\epsilon$, the $M_\theta$'s will still be global
slices, and the canonical line field will be
spanned by a vector field $\We$ of the form
$$\We = \ddth + W_1 ^{\epsilon}$$
where $W_1 ^{\epsilon}$ is tangent to the foliation
by the $M_{\theta}$. 
The $M_{\theta}$'s have a 1-parameter family
of contact structures, namely  
the intersections of the rank 3
distribution $(\D^{\epsilon}) ^2$ with their tangent
spaces.  

By the last lemma of the 
previous section, the time $t$ flow
 of the vector field $\We$
defines   contact diffeomorphisms
between slices $M_{\theta}$ and $M_{\theta +t}$.
Let us write $\phiet$ for this
map restricted to the bottom $M_0 = M$
of $\Omega$.  Thus
$\phiet:M_0 \to M_t$.
Using the unperturbed identification
(as the sections $V_{\theta}$)
of each of  the $M_t$ with $M$,  
we may think of $\phiet$ as a map
$M \to M$.  
For $t$ close to $0$ or
close to $\pi/2$  these maps $\phiet$ are
contact maps with respect to
the original contact structure.  This 
is   because the 
original  Engel structure has been
left undeformed near the top and bottom of
$\Omega$. It follows (see lemma 2) that  the
contactification of the perturbed Engel structure remains 
constant under perturbation and is
contact diffeomorphic to the original structure
$(M, \xi)$.

Let us describe this diffeomorphism explicitly. 
Let $\LLe$ denote the foliation by curves
defined by the perturbed
line-field.   Then we have  
that  
 $\Omega/ \LLe  \cong  M$ as contact manifolds. 
The identification is obtained by
mapping the leaf $\LLe (q)$ through a point $q$ to the
point  where that leaf   intersects the bottom $M_0 = M$,
 of $\Omega$. 
The projection
map $\pie: \Omega \to M$, 
whose fibers are the leaves of the perturbed foliation
$\LLe$,   can then be written 
as
$$\pie (m, \theta) = \phi^{\epsilon -1} _{\theta} (m).$$
And the induced contact structure, $\pie_* \De^2$,
equals the original contact structure.  
 
The intersection of the perturbed Engel field
$\De$ with the tangent space to the
$\Mth$ defines a line field, say $\elle$
which is canonically oriented.  This line
field is spanned by a vector field which together with
$\We$ spans $\De$.     The development map (see
\S 5.2 and 5.3) 
$$\Psie: (\Omega , \De) \to (S \xi  ,\D)$$
is now given by:
$$\Psie (q) = d \pie (q) (\elle (q)).$$
For $\theta$ close to $0$ or $\pi/2$
the line field $\elle$ on $\Mth$ is the one spanned
by our original $V(\theta)$. 

It follows from this discussion and
  the previous section
 that $\Psie$ is an Engel embedding
of the perturbed Engel manifold
with boundary, $(\Omega, \De)$
into the original Engel manifold
$S \xi $. Moreover the image
of $\Psie$ is $\Omega(V_0, V_1 ^{\epsilon})$
where
$$\Ve = (\phie_{\pi/2})^* V_1.$$ 

QED.  

\subsection{Engel automorphisms; Proof of the rest of Theorem 1}

In order to prove the assertions in the second paragraph
of theorem 1
we  need to know more about
Engel automorphisms.  
 Let $\phi: Q \to Q$ be an Engel
automorphism:  $\phi^* \D = \D$.  Then
$\phi$ preserves all the invariants of $\D$
so we have $\phi^* \D^2 = \D^2$
and $\phi^*L = L$.  
Let $U$ be any open set such that $U/{\cal L}$ is
a nice quotient.  For example, $U$ might be a small
tubular neighborhood of any local slice $M$ .  
Then $\phi$ induces a contact automorphism of the 
local contactifications, that is, a map from the local quotient 
$U/ {\cal L}$  to the local quotient 
$\phi(U)/{\cal L}$ 
It follows that the restriction of
$\phi$ to a local slice is
also a contact automorphism,  
$M \to \phi(M) = M^{\prime}$. 

\begin{lemma}  The restriction of 
an Engel isomorphism to a local slice
uniquely determines that isomorphism in
a neighborhood of the slice.  If
the slice is global then this restriction
uniquely determines the entire map.
\end{lemma}

{\sc proof}  Let $\phi$ and $\psi$ be two such isomorphisms,
and let $M$ denote  the slice.  Then $F= \psi ^{-1} \circ \phi$
is an automorphism of $X$ which is the identity
on $M$.  We will first show that $F$ is the identity in a 
neighborhood of $M$.    If $m \in M$ then, using the
construction in the BASIC EXAMPLE of \S 2
and lemma 2,  we can find 
canonical Engel coordinates $x,y,z,w$ such that
near $m$, the slice is defined by $w = 0$.
The remaining functions
$x,y,z$ are coordinates on $M$
so that $F$ has the form:
$F(x,y,z,w) = (x,y,z,g(x,y,z,w))$.
Now the Engel distribution is defined by
the Pfaffian system
$\theta_1 = \theta_2 = 0$ where
 $\theta_1 = dz - ydx$ and  $\theta_2 = dy -wdx$. 
For $F$ to be an Engel automorphism we must have
$F^* \theta_2 = 0  \hbox{ mod} \{\theta_1, \theta_2\}$, i.e.
$$dy-gdx = a \theta_1 + b \theta_2$$ for some functions
$a,b$.  Since no $dz$'s occur on the left hand side we
have $a = 0$.  Now write $g = w + h$. We obtain the equation
$ \theta_2 - hdx = b \theta_2$. Considering this equation
mod $\theta_2$,we obtain $h = 0$ so that $g = w$ and $F$
is the identity.

To prove the statement regarding global slices,
first use the fact that global slices are  local slices.
It follows that  the two maps agree in some neighborhood of the
slice.  
Now use the second paragraph of lemma 4 
which asserted that if two Engel automorphisms
agree on three points of such a leaf, then they
agree everywhere along that leaf. 

QED

\subsection{Proof of Theorem 2:} 
$M = M_0$ is a global slice for 
the Engel structure on $S \xi $.
It follows from  the lemma above and the
preceding discussions that 
every Engel isomorphism  
from the one domain to the other  is induced by
a contact map $ \psi: M \to M$. 
Such an isomorphism 
must take the boundary of the first
domain to the boundary of the second.
Being oriented, it must take the bottom boundary
to the boundary and top boundary to top boundary.
Hence it takes $V_0$ to $W_0$ and $V_1$ to $W_1$.

QED 
 
\subsection{Realization Lemma}
  
Here we prove the realization
lemma  (lemma 1) regarding contact isotopies.

{\sc Proof:}  Let $X = X(m;t) $ be a time-dependent
 vector field on $M$ which generates
the isotopy $\phi$ at time $t = \pi/2$.  We may insist 
that $X(t, m) = 0 $ for $t$ near zero and $t$ near 
$\pi/2$.  
The
original Engel distribution $\D$ is spanned by vector fields
$V(m, \theta) = cos(\theta) V_0 (m) + sin(\theta) V_1 (m)$ and
$\ddth$ where the contact structure $\xi$ on $M$ is 
spanned by the vector fields $V_0$ and $V_1$. 
 Define
a new  Engel structure 
$\E = span\{ V, W \}$ 
by declaring it to 
be framed by the same vector field $V(\theta)$, parallel
to the levels $\Mth$
and by the vector field
$$ W = \ddth + X (\theta, m)$$
transverse to the $\Mth$.  

Set 
$U(\theta, m) = {\partial \over {\partial \theta}} V$.
Then $U = -sin(\theta) V_0 (m) + cos(\theta) V_1 (m)$
and for each fixed $\theta$, the
fields $\{U(\theta, \cdot), V(\theta, \cdot) \}$
form a global frame for
the contact distribution.  One computes that 
$[W, V] = U+[X(\theta, \cdot), V(\theta, \cdot)]$
Since $X(\theta, \cdot)$ is an infinitesimal contact transformation
we can expand $[X, V] = f(\theta, m) V + g(\theta, m) U$
The condition that $\D^2$ has rank 3 is
precisely $g \ne -1$ anywhere.  If the perturbing field
$X$ is $C^1$-small enough then g will also be small and so   $g > -1$.  

In order to describe the perturbed $\D^2$,
we need some more information regarding the
contact generator $X$.
Let $\alpha$ be a choice
of contact one-form
(By a conformal change of $\alpha$
we may assume that $d \alpha(V, U) =1$.)
Let $Z$ be its
corresponding Reeb vector field, so
that $\alpha(Z) = 1, d \alpha(Z, \cdot ) = 0$.
Then $X$ must have the form
$h Z + X_h $ where
$h$ is a function and where $X_h$ is tangent to the
contact field and is determined by the
relation $dh|_{\xi} + i_{X_h} d \alpha |_{\xi} = 0$.
One computes that the perturbed $\D^2$,
which we will call $\E^2$,  is spanned by
$\ddth + h Z$ and $\xi$ or by $W$ and $\xi$
and hence has rank 3.  Since $[\xi, \xi]$ contains
$Z$ we see that $[\E, \E^2]$ has rank 4 so
that the distribution is indeed Engel.    
A straightforward calculation  
using the fact that $X$ is an infinitesimal
contact automorphism 
shows that 
$W$ spans the new characteristic line field:
$[W, \E ^2] \subset \E^2$. 
Let $\phi_s$ denote the time $s$ flow of
the time-dependent vector field $X$, starting at
time $0$. 
By construction, the time $s$ flow of $W$
is of the form $(m, 0) \mapsto (\phi_s (m), s)$,
when 
applied to a point $(m,0)$ at the `bottom level'
of $\Omega$  
 It follows that the bottom-to-top 
map is the desired contact map $\phi$.

QED

{\sc remark on smallness of the contact automorphism}

The smallness condition on $\phi$ is
solely a consequence of  the   constraint
$g \ne -1$ for its generator $X$.
This constraint is equivalent to $g > -1$
because $g$ is continuous and is zero  
 near the top and bottom of the domain. 
We attempt to give this last constraint
a geometric meaning. The initial choice of frame
$V_0, V_1$ determines an 
orientation  $A = V_0 \wedge V_1$ for $\xi$.
$V_0 \wedge V_1 = V \wedge U$. 
$V \wedge [X, V] = g V \wedge U$.
Then  $(V \wedge [X,V])/A = g$
is a measure of the amount of   `contact spinning'
induced by the infinitesimal contact
automorphism $X$.  The constraint $g > -1$
says  that the
contact transformation may not  
 ``spin'' the contact planes ``too much'' in 
the negative direction.  
This remark may be useful in the future in  
understanding large deformations 
or obstruction to existence of of 
Engel structures. (Eliashberg, private
communication.)

\subsection{Proof of Theorem 4 on
Zoll surfaces}
 
Theorem 4 follows directly from the following
lemma, combined with the realization lemma,
and the discussion following the statement
of the theorem. 

\ble
Let $g_t$ be a family of Zoll metrics
on the two-sphere. Let $V_t$ be the corresponding
Legendrian direction fields on
$STS^2$ defining
their geodesics.  Then there is
a contact diffeomorphism $\psi_t$
taking $V_t$ to $V_0$. 
\ele

The proof of this lemma in turn follows fairly
directly from part of the following theorem of Weinstein
\cite{Weinstein} (see also the appendix of \cite{Guillemin})
and the following lemma.   At first glance
this theorem of Weinstein appears
  closely related to ours.  However, 
as discussed after the statement of theorem 4, they are rather
different theorems.

\begin{backtheorem} [Weinstein] 
Let $g_t$ be a family of Zoll metrics
on the two-sphere. Let $Z_t$ be the corresponding
Hamiltonian vector fields defined on the
energy levels $\Sigma_t = \{ H_t = 1 \}$
for kinetic energy.  There exists
a one-parameter
family of contact maps $\beta_t: \Sigma_0 \to \Sigma_t$
such that $\beta_t ^* Z_t = Z_0$.
\end{backtheorem}

\ble Let $\xi_t$ be a deformation
of contact structures all of which
contain a fixed Legendrian
line field $L$.  Then
there is an isotopy $\phi_t$
with $\phi_t ^* \xi_t = \xi_0$
and $\phi_t ^* L= L$.
\ele

{\sc Proof of lemma 8.}

The metric $g_t$ induces a diffeomorphism
$F_t$ between $STS^2$ and $\Sigma_t$.
Namely, a tangent ray  $\ell = pos. span\{ {d x \over dt} \}$
attached at $x \in S^2$
is mapped to the unique  covector $p =  F_t (\ell)$
whose kernel $ker(p)$
defines a line   orthogonal to $\ell$ 
and normalized so that $p( \ell) >$
and so that the $g_t$ length of $p$ is 1.
It is important to observe that $F_t$ is
{\bf not} a contactomorphism.

We claim that the  pull-back $F_t ^* Z_t$ spans
the direction field $V_t$.  To see this let $(x(s), p(s))$
denote a solution to the geodesic equation written
in canonical variables and let $\dot x$ be the tangent
vector to the curve $x(s)$ in $S^2$.  It
suffices to show that $F_t (x, p) = (x, pos. span ( \dot x))$.
This follows directly from the   the Legendre transformation
relation
$\dot x^j = \Sigma g^{ji} p_i$.  For
the normal vector to a `hyperplane'
$\{ p_i = 0 \}$ is given by raising the indices of
$p_i$.  Consequently this normal vector equals the
tangent vector to the curve.  Moreover
the condition $H = \Sigma g_{ij}
p_i p_j = {1 \over 2}$ implies that $p_i \dot x^i = 1$
so that the orientations are as claimed.

Now use Weinstein's theorem to find a
diffeomorphism $\beta_t$ of $\Sigma_t$ such
that $\beta_t ^* Z_t = Z_0$.  (We do not
use the fact that $\beta_t$ is contact.)
Set $f_t = F_t ^{-1} \circ \beta_t \circ F_0$
and check that 
we have $f_t ^* V_t = V_0$. 

Set $\xi_t = f_t ^* \xi$, a deformation
of the contact structure $\xi$.  The
pair $(\xi_t, V_0)$ fulfills the hypothesis
of lemma 9.  Thus there
exists a contactomorphism
$\phi_t$ of $STS^2$ such that
$\phi_t^* \xi_t = \xi_0$
and $\phi_t ^* V_0 = V_0$.  The isotopy
$\psi_t = f_t \circ \phi_t$
now satisfies the conclusions of lemma 8.

QED

{\sc Proof of lemma 9.}

This can be proved using the standard 
deformation method as in Moser's proof
of  analogous results for volume
forms and symplectic
forms.   We look for an isotopy $\phi_t$ generated
by a time-dependent vector field $X_t$ which
we must solve for. Let $\theta_t$ denote contact forms for
$\xi_t$.   These forms are well-defined up
to a nonvanishing function $f_t$ so we must solve
$$\phi_t ^* (f_t \theta_t) =0.$$ 

Diferentiating this equation, using the assumption that
$X_t$ generates the flow, and factoring out
a common prefactor of $\phi_t ^*$ yields:
$$(L_{X_t} f_t) \theta_t + f_t( di_{X_t} \theta_t + i_{X_t} d \theta
_t) = - \dot f \theta_t - f_t \dot \theta_t.$$
Here dots denote $d /dt$ and we have
used Cartan's magic formula $L_X = di_X + i_X d$
for the Lie derivative of forms.   
We  view this
as a linear inhomogeneous equation
to be solved for $X_t, f_t$.  If they
can be solved then the flow $\phi_t$
of $X_t$ is the desired isotopy.  The only
assumption on the inhomogeneous terms
is that  $\dot \theta (L) = 0$.

To  proceed we assume that $X_t \subset \xi_t$.
Then $i_{X_t} \theta_t   =  0$
so that $di_{X_t} \theta_t = 0$.  We now break these
equations into their `horizontal' ($\xi_t$) and `vertical'
(Reeb) components by restricting both sides to $\xi_t$
and by applying both sides to the Reeb vector field.  The horizontal
equations read: $$i_{X_t} d \theta_t|_{\xi_t} = - \dot \theta_t
|_{\xi_t}$$ Since $d \theta_t|_{\xi_t}$ is symplectic
this equation has a unique solution $X_t \in \xi_t$.
Moreover the condition $\dot \theta (L) = 0$
implies that the solution $X_t$ {\bf must} lie in $L \subset \xi_t$.

Using this $X_t$ and the Reeb
field  we now write down the `vertical' equations.
 Recall that the
Reeb vector field $R_t$ is defined by the conditions $i_{R_t} \theta_t
= 1$ and $i_{R_t} d \theta_t = 0$.  Applying $i_{R_t}$
to 
both sides of our deformation equation
yields the equation
$$L_{X_t} f_t = -\dot f_t - f_t \dot \theta_t (R_t).$$
Set $g(t, x) = log(f_t (x)$.  The equation for f is
  satisfied if $g$ satisfies
$$\dot g = - L_{X_t} g - \dot \theta_t (R_t) .$$
This last equation
 is a linear inhomogeneous equation for $g_t (x) = g(t,x)$ which
can be solved by the method of variation of parameters.
The solution is 
$$g(t,x) = (\phi_{-t} ^* \int_0 ^t \phi_s ^* h(s, \cdot) ds ) (x),$$
where $h(t,x) =  \dot \theta_t (x) (R_t (x))$
and  where $\phi_t$ denotes the time t
flow of the non-autonomous vector field
$X_t$, starting at time 0.

We have found a solution $(X_t, f_t)$ for our system.
The flow $\phi_t$ yields the desired isotopy.  Observe
that $\phi_t ^* L= L$ since
$X_t$ is everywhere tangent to $L$.
QED 

\section{Appendix}

{\sc Geometry of  a pair of direction  fields}

Here we prove the lemma used in
our final section.  We restate the theorem
with our previous coordinates
$(x,y,p)$ replaced by $(x,z,y)$.  

\begin{lemma}  Let $\ell_0, \ell_1$ be two oriented
  line fields on a  three-manifold
 which  span a contact field.  Then in a
neighborhood of any point of   we can find coordinates
$x,y,z$ such that $$\ell_0 = span\{ \ddy \}$$
$$\ell_1 = span\{ \ddx + y \ddz + f \ddy\}$$
In these coordinates the contact field is
$ \{ dz - ydx = 0 \}$. If the line
fields depend continuously on a parameter
so do the coordinates and the
function $f$.  
\end{lemma}

{\sc remark.}The corresponding lemma for a single line
field $\ell_0$ contained in a contact
field is fairly well-known. 
It is the 3-dimensional case of the `Darboux's
theorem for Legendrian foliations'. See
p. 72 of \cite{ArnoldGivental}.  This
lemma can be found in some form in
the text on differential equations by
Arnol'd.   

{\sc proof:}  The proof proceeds
through a series of coordinate and frame changes.
The vector field spanning $\ell_0$ will
be denoted by $Y$ and the vector field spanning $\ell_1$
will be denoted by $X$.  The coordinates 
will be denoted $x,y,z$ and in each given step
we will go through one or more  coordinate changes:
$\bar x = \bar x (x,y,z), \bar y = \bar y (x,y,z),
\bar z = \bar z (x,y,z)$ and perhaps a frame
change:  $Y \to f Y, X \to g X$ with
$f, g$ positive functions.   

Step 1.  Find coordinates $x,y,z$ straightening
$\ell_0$, so that in this neighborhood
$\ell_0$ is spanned by $\ddy$.   Then $\ell_1$
is spanned by a vector field of the form
$f_1 \ddx + f_2 \ddy + f_3 \ddz$.  Without loss
of generality, we may assume $f_1 (0) > 0$.
(If $f_1 <0$ we can change the coordinate
$x$ to $-x$ to effect this change.  If $f_1 (0) = 0$
then it must be that $f_3 (0) \ne 0$ or otherwise
the two line-fields are not transverse.  In this
case we can switch the roles of $x$ and $y$.) 
Dividing by $f_1$ we have found coordinates
with
$$Y = \ddy$$
$$X = \ddx + f_2 \ddy + f_3 \ddz$$
By making a linear change of coordinates of the form
$$\bar x = x$$
$$\bar y = y + c_1 x$$
$$\bar z = z + c_2 x$$
we can cancel the constant term of $f_2$ and
$f_3$. This is done by taking $c_1 = -f_2 (0)$,
$c_2 = - f_3 (0)$, and using the chain rule.
Reverting to the unbarred coordinates, we can thus
impose the condition:
$$f_2(0) = f_3 (0) = 0$$

Step 2. We calculate that  
$$[Y,X] = {{\partial f_3} \over {\partial y}} \ddz
+ {{\partial f_2} \over {\partial y}} \ddy $$
The contact condition is that $X,Y, [X,Y]$ are linearly
independent everywhere, which means
that ${{\partial f_3} \over {\partial y}} \ne 0$.
We may assume that ${{\partial f_3} \over {\partial y}} > 0$
for if not, change the coordinate $z$ to $-z$ to insure
this.  Now define new coordinates
$$\bar x = x$$
$$\bar y = f_3(x,y,z) $$
$$\bar z = z $$
Then 
$$\ddx = \ddxbar + 
{{\partial f_3} \over {\partial x}}\ddybar$$
$$\ddy =  
{{\partial f_3} \over {\partial y}}\ddybar$$
 $$\ddz = \ddzbar + 
{{\partial f_3} \over {\partial z}}\ddybar$$
It follows that in the new coordinates,
$$Y = {{\partial f_3} \over {\partial y}}\ddybar$$
$$X = \ddxbar + \bar y \ddzbar + (f_3 + L_X f_2)
\ddybar.$$
Now, divide $Y$ by the positive function
${{\partial f_3} \over {\partial y}}$ and
set $f = (f_3 + L_X f_2)$ and
revert to unbarred coordinates.  We now have
$$Y = \ddy$$
$$X = \ddx + y \ddz + f \ddy$$
where $f$ is some arbitrary smooth
function of $x,y,z$. 

Finally, we observe that if the line
fields depend smoothly on a parameter,
then  all our coordinate changes
and our final function $f$  can
be made to depend smoothly on
this parameter as well.

\end{document}